\documentclass{amsart}
\usepackage{enumerate}
 \usepackage[latin1]{inputenc} 
\newtheorem{theorem}{Theorem}[section]

\theoremstyle{definition}

\theoremstyle{remark}
\newtheorem*{remark}{Remark}
\numberwithin{equation}{section}

\begin{document}

\title{A Question About Total Positivity and Newman's Fourier Transforms with Real Zeroes}

\author{Doug Pickrell}
\email{pickrell@math.arizona.edu}

\begin{abstract} Given a unitarily invariant ergodic measure on $\infty\times \infty$ Hermitian matrices, it is known that the characteristic function determines (and is determined by) a Polya frequency function $p(t)$. In turn the (finite) measure
$d\rho(u):=\frac{1}{p(-iu^2)}du$ has the property that the Fourier transform  $Z_b$ of $exp(-bu^2)d\rho(u)$ is an entire function and has real zeroes, for all $b\ge 0$; this is very close (but not identical) to a classification of such measures due to Newman. This raises the question of whether there is a direct connection between (e.g. the spectrum of) random Hermitian matrices and the reality of the zeroes of $Z_b$.

\end{abstract}
\maketitle

\setcounter{section}{-1}

\section{Introduction}

There are well-known conjectures (and piles of numerical evidence) relating the statistical properties of the eigenvalues of a large random Hermitian matrix (from the Gaussian Unitary Ensemble) and the zeroes of the zeta function along the critical line $Re(s)=\frac12$, see e.g. \cite{Mehta}. In this note we want to ask if there is a possible perturbation of this connection involving
Polya frequency functions (which correspond to totally positive kernels of the form $K(x,y)=f(x-y)$, in the terminology of \cite{Karlin}) and a class of functions characterized by Newman.

The rough outline is the following:

$$ \begin{matrix} f(a_{11})d\lambda(a_{11}) & \stackrel{IFT}{\rightarrow} & p(t) \\
\downarrow? & & \downarrow\\
Z_b(z) & \stackrel{FT}{\leftarrow} & d\rho(u)=\frac{1}{p(-iu^2)}d\lambda(u) \end{matrix} $$

As we will explain in the first section, $f$ is the pdf for the 1,1 entry of a random $\infty\times \infty$ Hermitian matrix, with a ergodic unitarily invariant probability distribution $\nu$ (The matrix coefficients of the random matrix make sense as random variables relative to $\nu$, but the spectrum does not a priori make sense as a random object). It turns out that $f$ is (generically)
a Polya frequency function. Schoenberg determined all the possible $f$, and these are expressed in terms of $p$, the characteristic function of $f$. Then, for reasons that are mysterious, we invert $p$ and replace $t$ with $-iu^2$. When we do this Schoenberg's class of functions transforms into a class of functions (or measures) studied by Newman (Newman's class of functions is slightly larger, but very close). Newman's measures are characterized by the fact that $Z_b$, the Fourier transform of $exp(-bu^2)d\rho(u)$, has real zeroes, for all $b\ge 0$.

The point of this note is simply to ask the following questions: (1) Is it possible to give a simple direct explanation for the fact that $f$ totally positive is essentially equivalent to $Z$ having real zeroes? 
(2) Is it possible to give an interesting interpretation to the zeroes of $Z$, e.g. are they related in some way to the eigenvalues of the original random matrix, distributed according to $\nu$? In particular does the diagram above explain why the zeroes are real? 

A basic example to bear in mind is the Gaussian unitary ensemble (GUE):

$$ \begin{matrix} (a_{11})_*GUE & \stackrel{IFT}{\rightarrow} & p(t)=e^{-t^2/2} \\
\downarrow? & & \downarrow\\
Z(z) & \stackrel{FT}{\leftarrow} & d\rho(u)=e^{-u^4/2}d\lambda(u) \end{matrix} $$
where GUE refers to the Gaussian (which is written heuristically as)
$$GUE=\frac1{\mathfrak z}exp(-\frac12 tr(A^2)) d\lambda(A)$$
What does the position of the real zeroes of $Z$ tell us about GUE?

\begin{remark}According to Maple, the Fourier transform of $exp(-u^{2m})$, for $m=2,...$, can be expressed in terms of hypergeometric functions. Maples does not seem able to find the Fourier transform of $exp(-u^4-bu^2)$.
The fact that the Fourier transform of $exp(-u^{2m})$ has real zeroes was proven by Polya; see
\cite{B} for perspective on Polya's work and references.  \end{remark}

\section{Ergodic Unitarily Invariant Measures on Hermitian Matrices and Total Positivity}

Let $e_1,e_2,...$ denote the standard orthonormal basis for $l^2$. Consider the direct limits
$Herm(\infty)=\lim Herm(n)$, realized concretely as
$$Herm(\mathbb Ce_1)\subset Herm(\mathbb Ce_1\oplus\mathbb Ce_2)\subset ...\subset Herm(\bigoplus_{1\le j\le n}\mathbb Ce_j)\subset...$$
and $U(\infty)=\lim U(n)$, realized as
$$U(\mathbb Ce_1)\subset U(\mathbb Ce_1\oplus\mathbb Ce_2)\subset ...\subset U(\bigoplus_{1\le j\le n}\mathbb Ce_j) \subset ...$$
There is a natural conjugation action
$$U(\infty)\times Herm(\infty)\to Herm(\infty):(g,X) \to gXg^{-1}$$
With respect to this action, any $U(\infty)$ invariant function of $X$ will be a function of the spectrum of $X$.

There is a natural identification of the dual of $Herm(\infty)$ with the linear space of
all Hermitian matrices $A=(a_{ij})_{1\le i,j<\infty}$ (without growth restriction), which we denote simply by $Herm$:
$$Herm(\infty)\times Herm\to \mathbb C: (X,A)\to trace(XA)$$
The action of $U(\infty)$ on the dual is identified with the natural action of $U(\infty)$ by conjugation on $Herm$,
$$U(\infty)\times Herm \to Herm:(g,A)\to gAg^{-1}$$
It does not make sense to talk about the spectrum of a general $A\in Herm$.

The set of $U(\infty)$ invariant probability measures on $Herm$, denoted $Prob(Herm)^{U(\infty)}$, is a convex set, and the extreme points are the ergodic measures. If $\nu\in Prob(Herm)^{U(\infty)}$, then its (inverse) Fourier transform
$$\hat{\nu}(X)=\int e^{itrace(XA)}d\nu(A)$$ is a function of the spectrum of $X$.

The following theorem is Proposition 5.9 in \cite{P}.

\begin{theorem} (a) $\nu\in Prob(Herm)^{U(\infty)}$ is ergodic if and only if the (inverse) Fourier transform is of the form
$$\hat{\nu}(X)=\prod_{j=1}^{\infty} p(t_j) $$
where the $t_j$ are eigenvalues of $X\in Herm(\infty)$.

(b) The possible $p$ are of the form
$$p(t)=e^{i\omega t}e^{-dt^2}\prod_{j=1}^{\infty}\left((1+id_jt)e^{-id_jt}\right)^{-1}$$
where $\omega\in\mathbb R$, $d\ge 0$, and $(d_j)\in l^2$.
\end{theorem}

Note that $p$ in part (b) has poles at the points $t=i/d_j$ along the imaginary real axis, for $d_j\ne 0$.

\begin{proof}Rough outline of (b): Write
$$p(t)=\int_{-\infty}^{\infty}e^{ita}f(a)da$$
Then (when $f$ is sufficiently smooth) it turns out (modulo some details) that $f$ satisfies (the totally positivity type condition
for Polya frequency functions)
$$det((f^{(i-1)}(\omega_j))_{1\le i,j\le n})\ge 0, \qquad \omega_1>\omega_2>...>\omega_n $$ for all $n$. A representation theorem
of Schoenberg (which can be found in the section of \cite{Karlin} on Polya frequency functions, although it is not easy to dig out) implies that $p$ has the desired form.

\end{proof}

Olshansky and Vershik \cite{OV} have given a different proof of the classification which does not use Schoenberg's representation theorem. As a consequence they are able to turn the result around to prove Schoenberg's theorem.
Bufetov \cite{Bufetov} has proven the surprising fact that a $U(\infty)$ invariant ergodic measure on $Herm$ is necessarily finite.

The basic example is $p(t)=exp(-t^2/2)$, which corresponds to GUE. The generic case can be viewed as a perturbation of
GUE. These perturbations are different from the perturbations related to quantum gravity which have been
prominent in the literature for the past 25 years. In quantum gravity one typically fixes large $N$ and considers a distribution on $Herm(N)$ of the form
$$\frac{1}{\mathcal Z}exp(-trace(g_2A^2+g_4A^4+..+A^{2m}))d\lambda(A)$$
This is $U(N)$ invariant, but (aside from GUE) this does not generally constitute a coherent family of measures as $N$ varies (In this framework, even for GUE, it is
necessary to let the coefficients depend on $N$ in an appropriate way to obtain asymptotic limits for the spectrum. I do not know if there is a known way to vary the coefficients with $N$ to obtain one of the measures in the theorem in the large $N$ limit - any limit will be unitarily invariant, but not necessarily ergodic).

\subsection{Newman's Classification}

Given $p$ as in the preceding subsection, we obtain a finite measure on $\mathbb R$
$$d\rho(u)=\frac{1}{p(-iu^2)}du $$
where
$$\frac{1}{p(-iu^2)}= e^{-\omega u^2-d u^4}\prod_{j=1}^{\infty}\left((1+d_j u^2)e^{-d_j u^2}\right) $$
It was known to Polya that
$$Z_b(z)=\int_{u=-\infty}^{\infty} exp(izu-bu^2)\frac{1}{p(iu^2)}du$$
has only real zeroes for all $b>0$. Newman \cite{Newman} proved a converse (a classification). As we vary $p$, these do not quite exhaust all the functions included in Newman's classification - one can multiply the function above by $u^{2m}$.

The fairy tale hope is that in some sense the zeroes of $Z$ are real because they are the (expected) spectrum of a Hermitian matrix, distributed according to $\nu$. How this might come about: Fix $N$. The projection of $\nu$ to $Herm(N)$, denoted $\nu_N$, is unitarily invariant and its Fourier transform is
$$\hat \nu_N(x)=\prod_{1}^N p(t_j)$$
where $t_j$ denotes the spectrum of $x\in Herm(N)$.

\section{Basic Examples}

Suppose that $p(t)=exp(i\omega t/2)$. In this case $\nu$ is a delta measure supported at a real multiple of the identity (these are fixed points for the conjugation action of $U(\infty)$. The spectrum is this multiple with infinite multiplicity.
$$d\rho(u)=exp(-\frac{\omega}{2} u^2)du$$
The partition function
$Z(z)$ is the Gaussian
$$Z(z)= \frac{1}{\sqrt{2\pi \omega}}exp(-\frac{1}{2\omega}z^2) $$
This does not have any zeroes. So this example seems kind of useless.

The following example somehow relates the Gaussian unitary ensemble and the Ising model with a magnetic field (In the latter context, the fact that we inserted an $i$ into Newman's definition possibly obscures the fact that $Z$ should be thought of as a partition function). What is special about the Ising model vs the other unitary minimal models? A square root which is quadratic.

Suppose that $p(t)=exp(-t^2/2)$. Then $d\rho(u)=exp(-u^4/2)$ and
$$Z_b(z)=\int_{u=-\infty}^{\infty} exp(izu-bu^2-u^4/2) du =exp(b^2/2)\int_{u=-\infty}^{\infty} exp(izu-\frac12 (u^2+b)^2 du$$
The completing the square does not appear to reduce the calculation of the integral to the $b=0$ case.

The function
$$Z_0(z)=\int_{u=-\infty}^{\infty} exp(izu-u^4/2) du$$
$$=\frac{1}{4\Gamma(3/4)}(-hypergeom([], [5/4, 3/2], (1/256)z^4)z^2\Gamma(3/4)^2+2\sqrt{2}\pi hypergeom([], [1/2, 3/4], (1/256)z^4))$$
is just a difference of hypergeometric functions.
It is an even function. But the graph along the real axis is crazy. For real relatively small $z$ it seems fairly reasonable. Then around $z=12$ it becomes extremely oscillatory. It is mainly positive, but it is occasionally zero, and occasionally huge. When it does have a zero, it appears to have lots of zeroes (because it is so oscillatory). Around $z=250$ it finally becomes small. It is very difficult to imagine what such a bizarre function could be telling us about GUE.

\subsection{Riemann xi function}

Riemann's xi function is defined by
$$\xi(z)=s(s-1)\pi^{-s/2}\Gamma(s/2)\zeta(s)/2, \qquad s=\frac12+iz$$
When $b=0$ this is the same as
$$\xi_b(u)=\int exp(izu-bu^2)F(t)dt$$
where $F$ is the even strictly even function
$$F(u)=\sum_{n=1}^{\infty}(4n^4\pi^2e^{9u/2}-6n^2\pi e^{5u/2})exp)-n^2\pi e^{2u}) $$
The fact that $F$ is even and positive is not clear - it is apparently a consequence of modular invariance of the Jacobi theta function,
see \cite{BPY}. There is a simpler presentation on Tao's blog, \cite{Tao}, in connection with his work on the de Bruijn-Newman constant.

$\xi$ is not equal to a $Z$. But does this framework suggest some interpretation for this function?


\begin{thebibliography}{99}

\bibitem{B} N. B. de Bruijn, The roots of trignometric integrals, Duke Math. J. 17 (1950) 197-226.

\bibitem{BPY} P. Biane, J. Pitman, and M. Yor, Probability laws related to the Jacobi theta function and
Riemann zeta functions and Brownian excursions, Bulletin (New Series) of the AMS, Vol. 38, Number 4 (2001) 435-465.

\bibitem{Bufetov} A. Bufetov, Finiteness of Ergodic Unitarily Invariant Measures on Spaces of Infinite Matrices,
ArXiv:1108.2737

\bibitem{Karlin}S. Karlin, Total Positivity, Stanford Univ. Press, Stanford, CA (1968).

\bibitem{Mehta}  M. L. Mehta, Random Matrices, 2nd edition (Academic, San Diego, 1991).

\bibitem{Newman}C. Newman, Fourier transforms with only real zeroes, Proc. AMS, Vol 61 No 2 (1976) 245-251.

\bibitem{OV} G. Olshanski and A. Vershik, Ergodic unitarily invariant measures on the space of infinite Hermitian matrices,
Contemporary Mathematical Physics. F. A. Berezin's memorial volume. Amer. Math. Transl. Ser. 2, vol. 175 (R. L. Dobrushin et al., eds), (1996) 137-175.

\bibitem{P}D. Pickrell, Mackey Analysis of infinite classical motion groups, Pac. J. Math., Vo. 150, No. 1 (1991) 139-166.

\bibitem{Tao} T. Tao, his blog post on de Bruijn-Newman constant.




\end{thebibliography}
\end{document}